\documentclass[11pt]{amsart}
\usepackage{amsmath,amsfonts,amsthm,amssymb}

\makeatletter

\newcommand{\theoremName}{Theorem}
\newcommand{\mtheoremName}{Main theorem}
\newcommand{\pmainname}{Proof of the main theorem}

\theoremstyle{plain}

\def\hm#1{#1\nobreak\discretionary{}{\hbox{\m@th$#1$}}{}}

\newtheorem {theorem}{\theoremName}
\newtheorem*{mtheorem}{\mtheoremName}
\newtheorem {lemma}{Lemma}
\newtheorem {corollary}{Corollary}
\newtheorem {definition}{Definition}

\newtheorem {remark}{Remark}
\newtheorem {example}{Example}
\newtheorem {proposition}{Proposition}

\let\@newpf\proof \let\proof\relax 
\newenvironment{proof}{\@newpf[\proofname]}{\qed\endtrivlist}

\DeclareMathOperator{\res}{Res}

\DeclareMathOperator{\Ind}{Ind}
\DeclareMathOperator{\Hom}{Hom}

\DeclareMathOperator{\Det}{Det}

\def\antiddots{\mathinner{\mkern1mu\raise\p@
    \vbox{\kern7\p@\hbox{.}}\mkern2mu
    \raise4\p@\hbox{.}\mkern2mu\raise7\p@\hbox{.}\mkern1mu}}
\makeatother

\DeclareMathAlphabet{\mathbbold}{U}{bbold}{m}{n}

\def\k{\mathbbold{k}}

\DeclareSymbolFont{rsfscript}{OMS}{rsfs}{m}{n}
\DeclareSymbolFontAlphabet{\mathrsfs}{rsfscript}

\DeclareFontFamily{OMS}{rsfs}{\skewchar\font'177}
\DeclareFontShape{OMS}{rsfs}{m}{n}{%
      <5> rsfs5
      <6> <7> rsfs7
      <8> <9> <10> rsfs10
      <10.95> <12> <14.4> <17.28> <20.74> <24.88> rsfs10
      }{}

\def\calA{\mathrsfs{A}}
\def\calB{\mathrsfs{B}}
\def\calC{\mathrsfs{C}}
\def\calD{\mathrsfs{D}}
\def\calE{\mathrsfs{E}}
\def\calF{\mathrsfs{F}}

\def\calO{\mathrsfs{O}}
\def\calP{\mathrsfs{P}}
\def\calQ{\mathrsfs{Q}}
\def\calR{\mathrsfs{R}}

\def\calV{\mathrsfs{V}}
\def\calW{\mathrsfs{W}}

\renewcommand{\L}{{\mathrsfs{L}\!ie}}
\renewcommand{\P}{\mathrsfs{P}}
\newcommand{\Com}{{\mathrsfs{C}\!om}}
\newcommand{\LL}{{\mathrsfs{L}\!ie_2}}
\newcommand{\CC}{{\mathrsfs{C}\!om_2}}
\newcommand{\PP}{{\mathrsfs{P}_2}}

\newcommand{\calHom}{\mathrsfs{E}\!nd}
\renewcommand{\ge}{\geqslant}
\let\phi\varphi

\title{
Character formulas\\
for the operad of two compatible brackets\\
and for the bihamiltonian operad}
\author{Vladimir Dotsenko, Anton Khoroshkin}
\date{}

\begin{document}
\begin{abstract}
We compute dimensions of the components for the operad of two compatible brackets and for the bihamiltonian operad. We also obtain character formulas for the representations of the symmetric groups and the $SL_2$ group in these spaces.
\end{abstract}
\maketitle

\section{Introduction}

\subsection{Description of the results.}
Let $V$ be a vector space equipped by two skew-symmetric operations such that each of them is a Lie bracket (i.e. satisfies the Jacobi identity), and moreover these two Lie brackets are compatible, i.e. any their linear combination is again a Lie bracket.
Such an algebraic structure is called an algebra with two compatible brackets.
As for an arbitrary algebraic structure, we can consider a free algebra of this type with $n$ generators $a_1,\ldots,a_n$. One of the natural problems concerning this algebra is to investigate how large it is. This question can be formulated more precisely as follows.
This algebra is graded (even multigraded), and graded components are finite-dimensional. What are the dimensions of the graded components? The common way to answer this question goes as follows.
First of all, study the $S_n$-module structure (w.r.t. the action induced by permuting the generators) on the ``multilinear part'' of the free algebra (i.e. all elements containing each of the generators exactly once). Then this information can be used to compute the dimensions of all graded components. In our case it is possible to prove that the dimension of the multilinear part of the free algebra is equal to $n^{n-1}$ and to compute the $S_n$-character for this space. In this paper we present these results.

We also derive similar results for ``bihamiltonian'' free algebras. A bihamiltonian algebra is an analogue of a Poisson algebra for the case of two compatible brackets. Namely, it is a vector
space equipped by two skew-symmetric mappings being two compatible brackets and a symmetric operation, being an associative commutative multiplication such that each of the brackets is its derivation. For a free algebra with $n$ generators of this type the dimension of the corresponding space is equal to
$(n+1)^{n-1}$, and it is also possible to compute the corresponding characters of the symmetric groups.

Conjecture on the dimension of the multilinear part of the free bihamiltonian algebra was stated by B.\,Feigin several years ago. At that moment it seemed that the corresponding spaces have much in common with the other representations of the symmetric groups of the same dimensions~--- the diagonal harmonics from the work~\cite{Ha}. The exact meaning of this phrase is not known yet; all evident ways to relate these spaces result in wrong statements.

\subsection{The machinery.}
Our main tool for computing dimensions and characters is the Koszul duality for the operads and the theory of Koszul operads developed by V.\,Ginzburg and M.\,Kapranov. It turns out that, though the operad of two compatible brackets is rather complicated, its Koszul dual is more visible. For any Koszul operad, it is possible to use the information on the dimensions of its components to obtain similar information for the dual operad:
\begin{proposition}[\cite{GK}]Let
$f_\calQ(x):=\sum_{n=1}^\infty\frac{\dim\calQ(n)}{n!}x^n$.
If the operad $\calQ$ is Koszul, then
$f_{\calQ}(-f_{\calQ^!}(-x))=x$.
\end{proposition}

\begin{example}
For the Lie operad $\L$ we have $\dim\L(n)=(n-1)!$, and so $f_\L(x)=-\ln(1-x)$.
For its Koszul dual $\Com$ we have $\dim\Com(n)=1$, so that $f_\Com(x)=\exp(x)-1$. It is well known that these operads are Koszul. The functional equation in this case reads as $-\ln(1+\exp(-x)-1)=x$.
\end{example}

A similar functional equation holds also for the generating functions of the symmetric group characters of the representations in the components of an operad. Thus the Koszulness of the operad of two compatible brackets turns out to be a very important ingredient of our work. We managed to prove it~--- in fact, there are even several ways to do it. First of them uses the ideas from the recent paper~\cite{V} establishing some connection between the Koszul operads and the Cohen--Macaulay partially ordered sets. Another way uses the criterion for the Koszulness in terms of distributive lattices (similar to the criterion for quadratic algebras) found by the second author.
(This proof is discussed in details in~\cite{Khor}.) After the proof of the Koszulness, all the calculations are just power series inversions.

\subsection{Plan of the paper.}
Throughout the paper we assume that the reader is familiar with the main operadic notions. Still we remind briefly some of them when they appear in the text.

In Section~\ref{LL-def} we define our operads and establish their basic properties. Here we avoid general operadic definitions, appealing to the intuitive understanding of the operads.
In Section~\ref{Koszul} we remind necessary notions of the Koszul duality for operads and define the generating functions for the characters. Calculating these generating functions for our operads is the main result of this paper.
In Section~\ref{poset-proof} we prove the Koszulness of our operads.
In Section~\ref{chi-calc} we compute the generating functions for the characters of our operads using functional equations on these generating functions.
In Section~\ref{monomial} we discuss the monomial basis in the multilinear part of the free algebra with two compatible brackets. The dimension formula helps us to prove that the family of monomials introduced by M.\,Bershtein~\cite{Ber}\footnote{As we were told by P.\,Etingof, another family of monomials which turns out to be a basis was independently discovered by Fu Liu.} is a basis.
In Section~\ref{mult} we list the current results on the decomposition of the components of our operads into sum of irreducibles. 
Finally, Section~\ref{basic} is auxiliary. It contains the definitions related to power series and residues and one simple calculation.

The ground field for all vector spaces and algebras throughout this paper is an arbitrary field 
$\k$ of zero characteristic. Sometimes we suppose that $\k$ is algebraically closed; we will mention it explicitly then.

\subsection{Acknowledgements.}
We are grateful to our teacher B.\,Feigin for stating the problem, useful discussions and significant help in editing this paper. We are also indebted to A.\,N.\,Kirillov who told us about the paper~\cite{V}, and to V.\,Ginzburg and P.\,Etingof for the discussion of our results.

The work of the first author is supported by the grant of the President of Russian Federation
2044.2003.2 and the INTAS grant 03-3350. The work of the second author is supported by the RFBR grant 04-01-00637 and the INTAS grant 03-3350. 

\section{Definitions of the operads.}\label{LL-def}

\subsection{Operads $\Com$, $\L$ and $\P$.}
We remind several standard definitions. The operad $\L$ is generated by a skew-symmetric binary operation $\{\cdot,\cdot\}$ with one quadratic relation: we want the Jacobi identity~\eqref{Jacobi} to be satisfied in each algebra over this operad. Thus an algebra over this operad is a Lie algebra.
The operad $\Com$ is generated by a symmetric binary operation $\star$ with one quadratic relation: the associativity law for this operation. An algebra over this operad is an associative commutative algebra.

The Poisson operad $\P$ is a little bit more complicated. It is generated by a symmetric operation $\star$ and a skew-symmetric operation $\{\cdot,\cdot\}$; the symmetric operation generates a suboperad of $\P$ isomorphic to $\Com$, the skew-symmetric operation generates a suboperad isomorphic
to $\L$ and the relations between these operations mean that the skew-symmetric one is a derivation of 
the symmetric one (``the Leibniz rule for differentiating a product''):
 $$
\{a,b\star c\}=\{a,b\}\star c+b\star\{a,c\}.
 $$
A algebra over this operad is called a Poisson algebra; an example of such an algebra is the algebra of functions on a manifold with the ordinary product and a bracket defined by a Poisson bivector field.

\subsection{Operads $\LL$ and $\CC$.}
The operad $\LL$ (called also the operad of two compatible brackets) is generated by two skew-symmetric operations (brackets) $\{\cdot,\cdot\}_1$ and $\{\cdot,\cdot\}_2$. The relations in this operad mean that any linear combination of these brackets satisfies the Jacobi identity.
It is equivalent to the following identities in each algebra over this operad: the Jacobi identity for each of the brackets
\begin{gather}\label{Jacobi}
\{a,\{b,c\}_1\}_1+\{b,\{c,a\}_1\}_1+\{c,\{a,b\}_1\}_1=0,\\
\{a,\{b,c\}_2\}_2+\{b,\{c,a\}_2\}_2+\{c,\{a,b\}_2\}_2=0,
\end{gather}
and a six-term relation between the brackets
\begin{multline}\label{6term}
\{a,\{b,c\}_1\}_2+\{b,\{c,a\}_1\}_2+\{c,\{a,b\}_1\}_2+\\
+\{a,\{b,c\}_2\}_1+\{b,\{c,a\}_2\}_1+\{c,\{a,b\}_2\}_1=0.
\end{multline}

The operad of two commutative products $\CC$ is generated by two symmetric binary operations (products) $\star_1$ and $\star_2$ such that in any algebra over this operad we have the associativity conditions
\begin{gather}
a\star_1 (b\star_1 c)=(a\star_1 b)\star_1 c,\\
a\star_2(b\star_2 c)=(a\star_2 b)\star_2 c,
\end{gather}
and also the identities
\begin{multline}
a\star_1(b\star_2 c)=a\star_2(b\star_1 c)=b\star_1(a\star_2 c)=\\
=b\star_2(a\star_1 c)=c\star_1(a\star_2 b)=c\star_2(a\star_1 b)
\end{multline}
meaning that the six elements obtained from the given three by taking two distinct products are pairwise equal.

Note that the components of these operads are equipped by an $SL_2$-action (besides the action of the symmetric groups), since $SL_2$ acts on binary operations that generate these operads.

\begin{proposition}\label{CCdim}
The space $\CC(n)$ is $n$-dimensional.
As a representation of $S_n$ it is isomorphic to the direct sum of $n$ copies of a trivial representation. As a representation of $SL_2$ it is isomorphic to the $n$-dimensional irreducible representation $L(n)$.
\end{proposition}

\begin{proof}
Notice that the dimension of $\CC(n)$ is at most~$n$: for each
$k=0$, $1$, \ldots, $n-1$ all monomials obtained by $k$ products of the first type and $n-1-k$ products of the second type are pairwise equal. On the other hand, this suggests an explicit realization of the free algebra over this operad. Namely, an element of the basis is indexed by a multiset whose elements are the generators of the algebra and a nonnegative integer less than the cardinality of the multiset (the number of the products of the first type in the monomial); the definition of the operations is clear. It follows that the dimension is equal to~$n$. The statements concerning representations are straightforward (for $SL_2$ the simplest way is just to compute the weights which turn out to be $n-1$, $n-3$, \ldots, $1-n$, and this is exactly what we need).
\end{proof}

\subsection{The operad $\PP$.}
The bihamiltonian operad $\PP$ is generated by three operations: two skew-symmetric ones
($\{\cdot,\cdot\}_1$ and $\{\cdot,\cdot\}_2$) and a symmetric one ($\star$). The relations
on these operations mean that the skew-commutative operations are two compatible Lie brackets
(that is, relations \eqref{Jacobi}--\eqref{6term} hold), the commutative operation is an associative product and each of the brackets is its derivation:
\begin{gather}
\{a,b\star c\}_1=\{a,b\}_1\star c+b\star\{a,c\}_1,\label{Leibniz}\\
\{a,b\star c\}_2=\{a,b\}_2\star c+b\star\{a,c\}_2.
\end{gather}

\section{Operads: a summary.}\label{Koszul}

\subsection{$\mathbb{S}$-modules and operads.}
By $\mathbb{S}\mathrm{-mod}$ we denote the category of $\mathbb{S}$-modules, i.e. collections of vector spaces $\{\calV(n), n\ge1\}$, where each $\calV(k)$ is an $S_k$-module. The set of all morphisms between  $\calV$ and $\calW$ in this category is, by definition,
 $$
\Hom_{\mathbb{S}\mathrm{-mod}}(\calV,\calW)=\bigoplus\limits_{n\ge1}\Hom_{S_n}(\calV(n),\calW(n)).
 $$
Define the direct sum, the tensor product and the dual module by the formulas $(\calV\oplus\calW)(n)\hm=\calV(n)\oplus\calW(n)$,
$(\calV\otimes\calW)(n)\hm=\calV(n)\otimes\calW(n)$, $(\calV^*)(n)=(\calV(n))^*$.

An important example of an $\mathbb{S}$-module is the module $\Det$, for which $\Det(n)$ is a sign representation of~$S_n$. We use this module to define a version of a dual module which we actually use: $\calV^\vee=\calV^*\otimes\Det$; it is the ordinary dual twisted by the sign representation.
In some cases we consider the differential graded  $\mathbb{S}$-modules; all the above constructions for them are defined analogously. The graded analogue of $\Det$ is denoted by $\calE$;
the space $\calE(n)_{1-n}$ is one-dimensional and is the sign representation of the symmetric group while all other spaces $\calE(n)_k$ are equal to zero.

Each $\mathbb{S}$-module $\calV$ gives rise to a functor from the category $\calF\!in$ of finite sets (with bijections as morphisms) to the category of vector spaces. Namely, for a set $I$ of cardinality $n$ let
 $$
\calV(I)=\k\Hom_{\calF\!in}([n],I)\otimes_{\k S_n}\calV(n).
 $$
(Here $[n]$ denotes the ``standard'' set $\{1,2,\ldots,n\}$.)

For $\mathbb{S}$-modules $\calV$ and $\calW$ define the composition $\calV\circ\calW$ as
\begin{equation}
(\calV\circ\calW)(n)=\bigoplus\limits_{m=1}^n\calV(m)\otimes_{\k S_m}
\left(\bigoplus\limits_{f\colon [n]\twoheadrightarrow [m]}\bigotimes\limits_{l=1}^m\calW(f^{-1}(l))\right),
\end{equation}
where the sum is taken over all surjections $f$. The following definition is more visible but less invariant:
\begin{equation}
(\calV\circ\calW)(n)=\bigoplus\limits_{m=1}^n
\calV(m)\otimes_{\k S_m}\left(\bigoplus\limits_{i_1+\ldots+i_m=n}\bigotimes\limits_{l=1}^m\calW(i_l)\right).
\end{equation}

This operation endows the caterory of $\mathbb{S}$-modules with a structure of a monoidal category. An operad is a monoid in this category. See~\cite{MSS} for a more detailed definition.

A set-theoretic operad can be defined similarly. The only difference is in replacing $\mathbb{S}$-modules by $\mathbb{S}$-sets (i.e. collections of sets $\{M_k, k\ge1\}$ with the symmetric groups action; the tensor product should be replaced by the Cartesian product).

Given a set-theoretic operad, one can construct an ordinary operad, replacing the sets by their linear spans and extending the composition maps to linear combinations by linearity.

Let $V$ be a vector space. By the definition, the operad of linear mappings $\calHom_V$ is a collection $\{\calHom_V(n)=\Hom(V^{\otimes n},V), n\ge1\}$ of all multilinear mappings from $V$ to itself with the obviuos composition maps.

Using the operad of linear mappings, we can define an algebra over an operad $\calO$;
a structure of such an algebra on a vector space is a morphism of the operad $\calO$ to the corresponding operad of linear mappings. Thus an algebra over an operad $\calO$ is a vector space $W$ together with a collection of mappings $\calO(n)\otimes_{\k S_n} W^{\otimes n}\to W$ with obvious compatibility conditions. A free algebra generated by a vector space $X$ over an operad $\calO$ is (isomorphic to)
$\bigoplus_{k=1}^\infty\calO(n)\otimes_{\k S_n}X^{\otimes n}$.

To simplify the definitions, we consider here only operads~$\calQ$ with $\calQ(1)=\k$.

\subsection{Operads defined by generators and relations.}
A free operad $\calF_\calO$ generated by an $\mathbb{S}$-module $\calO$ (with $\calO(1)=0$) is defined as follows. A basis in this operad consists of some species of trees. These trees have a distinguished root (of degree one). A tree belonging to $\calF_\calO(n)$ has exactly $n$ leaves.
Finally, this tree has internal vertices (neither leaves nor the root) labeled by the basis of the set of generators $\calO$, and a vertex with $k$ siblings is labeled by an element from $\calO(k)$.
The unique tree whose set of internal vertices is empty generates a one-dimensional space  $\calF_\calO(1)$. The composition of a tree $t$ with $l$ leaves and the trees
$t_1,\ldots,t_l$ glues the roots of $t_1,\ldots,t_l$ to the corresponding leaves of $t$ (so that two glued egdes become one edge, and the corresponding vertex becomes an internal point of this edge).

Free operads are used to define operads by generators and relations. Let $\calO$ be an $\mathbb{S}$-module, $\calR$ be an $\mathbb{S}$-submodule in $\calF_\calO$.
An (operadic) ideal generated by $\calR$ in $\calF_\calO$ is the linear span of all trees such that at least one internal vertex is labeled by an element of~$\calR$.
An operad with generators $\calO$ and relations~$\calR$ is the quotient of the free operad modulo this ideal.

\subsection{Koszul duality for operads.}

Let an operad $\calQ$ be defined by a set of binary operations $\calO$ with 
quadratic relations $\calR$ (that is, they involve ternary operations obtained by composition maps from the given binary ones). In this case $\calQ$ is called quadratic. For quadratic operads, there is an analogue of the Koszul duality for quadratic algebras. This duality assigns to a quadratic operad $\calQ$ an operad $\calQ^!$ with generators $\calO^\vee$ and the annihilator of $\calR$ under the natural pairing as the space of relations. The property $(\calQ^!)^!\simeq\calQ$ holds here as well.

\begin{proposition}
$\Com^!\simeq\L$, $\P^!\simeq\P$, $\LL^!\simeq\CC$.
\end{proposition}

\begin{proof}
First two isomorphisms are well known (being established, for example, in~\cite{GK}), and the last one is similar to the first one.
\end{proof}

The cobar complex $\mathbf{C}(\calQ)$ for an operad $\calQ$ is a free operad with generators
$\{\calQ^*(n), n\ge2\}$ equipped with a differential $d$ with $d^2=0$. We will give here a sketch of the definition, referring the reader to \cite{MSS} for the details (including the grading on the free operad and the appropriate choice of signs in the formula for $d$).

Notice that from the operad structure on $\calQ$ we can obtain the maps
 $$\phi^*_{m,n,k}\colon\calQ^*(m+n-1)\to\calQ^*(m)\otimes\calQ^*(n)$$
dual to the composition maps
\begin{multline}
\phi_{m,n,k}\colon\calQ(m)\otimes\calQ(n)\simeq\\
\simeq\calQ(m)\otimes\calQ(1)^{\otimes (k-1)}\otimes\calQ(n)\otimes\calQ(1)^{\otimes (m-k)}\to\calQ(m+n-1)
\end{multline}
in $\calQ$. A differential on the generators of the free algebra is defined as a certain sum of these maps with alternating signs, and it can be extended by an operadic analogue of the Leibniz identity (with some signs that take into account the grading) to compositions of generators.

Once again we use twisting by the sign, now to get another version of the cobar complex:
$\mathbf{D}(\calQ)=\mathbf{C}(\calQ)\otimes\calE$.
The zeroth cohomology of~$\mathbf{D}(\calQ)$ is isomorphic to the operad~$\calQ^!$.

\begin{definition}
An operad $\calQ$ is Koszul if $H^i(\mathbf{D}(\calQ))=0$ for~$i\ne0$.
\end{definition}
\begin{remark}
It is well known (\cite{GK}, \cite{M}) that operads $\L$, $\Com$ and $\P$ are Koszul.
\end{remark}

Recall that for a quadratic operad~$\calQ$ the $\mathbb{S}$-module
$\calQ\circ ((\calQ^!)^\vee)$ can be equipped by a differential which endows it with a structure of a cochain complex (see \cite{MSS}).
This complex is called the Koszul complex for an operad $\calQ$ (and is similar to the Koszul complex for a quadratic algebra).
\begin{proposition}[\cite{GK}]\label{equiv}
Let $\calQ$ be a quadratic operad. The following conditions are equivalent:
\begin{itemize}
\item $\calQ$ is Koszul.
\item $\calQ^!$ is Koszul.
\item For $n>1$ the subcomplex $\calQ\circ ((\calQ^!)^\vee)(n)$ of the Koszul complex for~$\calQ$ is acyclic.
\end{itemize}
\end{proposition}

\subsection{Distributive laws.}
Here we briefly remind the main result of the Markl's work~\cite{M}.
Let $\calA$ and $\calB$ be two quadratic operads. Denote by $\calO_\calA$, $\calO_\calB$ their generators and by $\calR_\calA$, $\calR_\calB$ the spaces of relations. Denote by
$\calO_\calA\bullet\calO_\calB$ a subspace in the free operad $\calF_{\calO_\calA\oplus\calO_\calB}$ spanned by the elements $\phi(1,\psi)$ where $\phi\in\calO_\calA$, $\psi\in\calO_\calB$ (and $1$ stands for the identity unary map corresponding to the tree without internal vertices). The notation $\calO_\calB\bullet\calO_\calA$ has the same meaning. Suppose that we have a map
 $$
d\colon\calO_\calB\bullet\calO_\calA\rightarrow\calO_\calA\bullet\calO_\calB.
 $$
Consider an operad $\calC$ with generators $\calO_\calC=\calO_\calA\oplus\calO_\calB$ and relations
$\calR_\calC\hm=\calR_\calA\oplus\calD\oplus\calR_\calB$, where $\calD=\{x-d(x)\mid x\in\calO_\calB\bullet\calO_\calA\}$. It is easy to see that the natural inclusion 
 $$
\calF_{\calO_\calA}\circ\calF_{\calO_\calB}\hookrightarrow\calF_{\calO_\calC}
 $$
gives rise to a morphism of $\mathbb{S}$-modules $\xi\colon\calA\circ\calB\rightarrow\calC$. 
Let the degree of elements from $\calO_\calA$ be equal to $(1,0)$, and the degree of elements from $\calO_\calB$ be equal to $(0,1)$. Thus our $\mathbb{S}$-modules become bigraded, and the morphism respects the grading.

The next proposition is proved in~\cite{M}.

\begin{proposition}\label{distrib}
Suppose that $\calA$ and $\calB$ are Koszul, and $\xi$ is an isomorphism on the homogeneous components of degrees $(2,1)$ and $(1,2)$. 
Then
\begin{itemize}
\item[(i)] $\xi$ is an isomorphism of bigraded $\mathbb{S}$-modules $\calA\circ\calB$
and $\calC$.
\item[(ii)] The operad $\calC$ is Koszul.
\end{itemize}
\end{proposition}

If the condition from the proposition holds, the map $d$ is called a distributive law between $\calO_\calB\bullet\calO_\calA$ and $\calO_\calA\bullet\calO_\calB$.

\subsection{Generating functions and characters.}
As we have already mentioned in the introductory part, to each operad (and more generally, to each $\mathbb{S}$-module) $\calQ$ one can assign a formal power series~--- the exponential generating function for the dimensions
 $$
f_\calQ(x)=\sum_{n=1}^\infty\frac{\dim\calQ(n)}{n!}x^n,
 $$
such that for a Koszul operad 
 $$
f_{\calQ}(-f_{\calQ^!}(-x))=x.
 $$
This functional equation is an immediate corollary of a functional equation relating more general generating functions which we will define here.

Character of a representation $M$ of the symmetric group $S_n$ can be identified (\cite{Mac})
with a symmetric polynomial $F_M(x_1,x_2,\ldots)$ of degree~$n$ in infinitely many variables.
To each $\mathbb{S}$-module $\calV$ we assign an element $F_\calV(x_1,\ldots,x_k,\ldots)\hm=
\sum_{n\ge1}F_{\calV(n)}(x_1,\ldots,x_k,\ldots)$ of the ring of symmetric functions~$\Lambda$. This ring is the completion of the ring of symmetric polynomials in infinitely many variables with respect to the valuation defined by the degree of a polynomial. It is isomorphic to the ring of formal power series in (infinitely many) variables of degrees $1,2,3,\ldots$ (for which one can take, for example, the Newton power sums $p_1,\ldots,p_n,\ldots$). The series $F_\calV$ is a generating series for the symmetric groups' characters. Namely, multiplying the coefficient of $p_1^{n_1}\ldots p_k^{n_k}$ by
$1^{n_1}n_1!\ldots k^{n_k}n_k!$, we obtain the value of the character of $\calV(n)$ on a permutation containing $n_1$ cycles of length~$1$, \ldots, $n_k$ cycles of length~$k$ in the decomposition into disjoint cycles.

Note that $f_\calV(x)=F_\calV|_{p_1=x,p_2=p_3=\ldots=0}$.

If $\calV$ is equipped by the action of a group $G$ which commutes with the symmetric groups, then for each $n$ the space $\calV(n)$ is a representation of $S_n\times G$. In this case we assign to $\calV$ an element of the ring $\Lambda_{G}$ of symmetric functions with values in the character ring of $G$ (or, in other words, a character of $G$ with values in symmetric functions). We denote this element by $F_\calV(x_1,\ldots,x_n,\ldots;g)$ with $g\in G$.

Extend the definition of $F_\calV$ on differential graded modules. Namely, for such a module $\calV=\oplus_i\calV_i$ let $F_\calV=\sum_i(-1)^iF_{\calV_i}$ (the Euler characteristic of $\calV$).

\subsection{Functional equation for characters.}
Let us remind the definition of the plethysm for the symmetric functions.
\begin{definition}
Fix a symmetric function $H(x_1,x_2,\ldots)$. Plethysm (plethystic substitution of~$H$) is a
$\k$-linear ring homomorphism $\Lambda\to\Lambda$ (the image of $F$ under this homomorphism is denoted by $F\circ H$) defined by $p_n\circ H\hm=H(x_1^n,x_2^n,\ldots)$.
\end{definition}

\begin{proposition}[\cite{GK}]
For each two $\mathbb{S}$-modules $\calV$, $\calW$ we have
 $$
F_{\calV\circ\calW}=F_{\calV}\circ F_{\calW}.
 $$
\end{proposition}

Extend this standard defintion on elements of $\Lambda_{G}$ in a following evident way.
\begin{definition}
Fix $H(x_1,x_2,\ldots;g)\in\Lambda_G$. Plethysm (plethystic substitution of~$H$) is a ring homomorphism $F\mapsto F\circ H$ of $\Lambda_G$ to itself that is linear over the character ring of $G$ and is defined on symmetric functions by
$p_n\circ H=H(x_1^n,x_2^n,\ldots; g^n)$.
\end{definition}
The formula for the character of the composition can be proved analogously.

From the definition it is clear that
 $$
p_n\circ (H(p_1,p_2,\ldots,p_k,\ldots;g))=H(p_n,p_{2n},\ldots,p_{kn},\ldots;g^n).
 $$

Define by $\varepsilon$ the linear over the character ring of $G$ involution of $\Lambda_G$ mapping  $p_n$ to $-p_n$. 

\begin{theorem}\label{FunEq}
Suppose that the operad $\calQ$ is Koszul. Then the following equality holds in $\Lambda_G$: \begin{equation}
\varepsilon(F_{\calQ})\circ\varepsilon(F_{\calQ^!})=p_1.\\
\end{equation}
\end{theorem}
The proof is similar to the proof of the corresponding statement in \cite{GK};
to prove this formula, one compares two different calculations of the Euler characteristic of the Koszul complex for $\calQ$.

\section{Koszulness and Cohen--Macaulayness.}\label{poset-proof}

In this section we prove that our operads are Koszul. The proof uses the results relating Koszul operads and Cohen--Macaulay partially ordered sets.
\begin{definition}
A partially ordered set\footnote{Throughout this paper we use the standard abbreviation ``poset''.} is a set $M$ equipped by a transitive binary relation $<$ such that no pair of elements $a,b\in M$ satisfies both $a<b$ and $b<a$. This relation is called the order relation on $M$.
\end{definition}
If for the elements $a,b$ we have $a<b$ we say that $a$ is less than $b$, and $b$ is greater than $a$.
\begin{definition}
\begin{enumerate}
\item A chain in a poset $M$ is an increasing sequence of elements $(a_0<\ldots<a_k)$.
\item Let $x<y$ be two elements of a poset $M$. The interval $(x,y)$ is a subset of $M$ consisting of all $z$ such that $x<z<y$. The segment $[x,y]$ is a subset $(x,y)\cup\{x\}\cup\{y\}$.
\item If $x<y$ and $(x,y)=\varnothing$ we say that $y$ covers $x$ and denote it by $x\prec y$.
\end{enumerate}
\end{definition}
In this section we consider only set-theoretic operads and linear spans of such operads. We use the same notation for an operad and its linear span if it is clear from the context which of them is considered. Fix a convenient way to work with an operad and its linear span simultaneously. Namely, for a realization of the $n$th set of the operad $\calA$ we take the set of multilinear monomials in the free algebra with $n$ generators over the linear span of this operad.

Our main example is the operad $\CC$. The basis in the multilinear part of the free algebra is due to the relations a quotient of a finite set modulo some equivalent relation. Hence this operad is a linear span of a set-theoretic operad.

Let $\calA$ be a set-theoretic operad. Let us construct from this operad a collection of posets
$\Pi_n(\calA)$ as follows. An element of $\Pi_n(\calA)$ is a pair $(I,\psi)$, where
$I$ is a partition $I_1\sqcup I_2\sqcup\ldots\sqcup I_k$ of an $n$-element set into a disjoint union of subsets and $\psi$ is a collection of elements belonging to the operad: for each $s=1,\ldots,k$ we have $\psi_s\in\calA(I_s)$. By the definition, $(I,\psi)$ is less than $(J,\phi)$ if $I$ is obtained from $J$ by refinement (partitioning some of the subsets; $J_s=\cup_{l=1}^{m_s} I_{a(l,s)}$), and the ``operations'' $\phi$ are obtained from $\psi$ by taking the operadic compositions with some elements
$\eta_s\in\calA(m_s)$: $\phi_s=\eta_s(\psi_{a(1,s)},\ldots,\psi_{a(m_s,s)})$.

Koszulness of the linear span of this operad is related to some homological properties of the constructed posets.
\begin{definition}
The order complex of a poset $M$ is a vector space $C_*(M)$ with basis indexed by chains
$(m_0<\ldots<m_k)$ and the differential $d$:
 $$
d(a_0<\ldots<a_k)=\sum_i(-1)^i(a_0<\ldots<a_{i-1}<a_{i+1}<\ldots<a_k).
 $$
\end{definition}
\begin{definition}
The homology of the order complex is called the homology of the poset $M$ and is denoted by $H_*(M)$.
\end{definition}
\begin{definition}
A poset $M$ is called Cohen--Macaulay if $H_i(M)=0$ for all positive $i$ less than the length of maximal chains in $M$.
\end{definition}
The following result appears in \cite{V}.
\begin{proposition}\label{poset-operad}
Suppose that for any elements $\alpha_i\in\calA(m_i)$, $i=1,\ldots,k$ of a set-theoretic operad  $\calA$
\begin{itemize}
\item[(*)]
the map $\calA(k)\to\calA(m_1+\ldots+m_k)$,
$\beta\mapsto\beta(\alpha_1,\ldots,\alpha_k)$, is injective.
\end{itemize}
Then the linear span of this operad is Koszul if and only if all posets $\Pi_n(\calA)$ are Cohen--Macaulay.
\end{proposition}
This proposition was successfully applied in the paper \cite{CV} containing among other results the proof for the Cohen--Macaulayness of the poset of pointed partitions of a finite set.
\footnote{These finite posets are the posets assigned to the set-theoretic operad $\calP erm$ (see\,\cite{V}). In general, the operads $\calP erm$ and $\CC$ have much in common; for example, they have the same generating series for the dimensions (but the symmetric groups actions are different, and it leads to the crucial difference between the dual operads).} The idea of the proof for our case has much in common with that of \cite{CV}.

The following proposition is evident.
\begin{proposition}
The condition~(*) holds for the operad~$\CC$.
\end{proposition}
\begin{definition}
Call a poset $M$ upper semimodular if for all $a,b,c\in M$ such that $c\prec a,b$ there exists $d\in M$ such that $a,b\prec d$.\par
A poset $M$ is called totally upper semimodular if any its segment is upper semimodular.
\end{definition}
The results of paper \cite{BW} lead (see\,\cite{CV}) to a following criterion of Cohen--Macaulayness.
\begin{proposition}
A totally upper semimodular poset $M$ is Cohen--Macaulay.
\end{proposition}

Consider the sets $\Pi_n=\Pi_n(\CC)$ more precisely. An element of $\Pi_n$ is a set of monomials in a free algebra with $n$ generators over $\CC$ where each generator appears (in total) once. The refinement corresponds to replacing some monomials from the set by their factors.
\begin{theorem}
For each $n$ the poset $\Pi_n$ is Cohen--Macaulay.
\end{theorem}
\begin{proof}
Let us prove that $\Pi_n$ is totally upper semimodular. We start with the maximal (w.r.t. the length) segments $\Pi_{n,s}$ of this set. The maximal element of $\Pi_{n,s}$ is the monomial using $s$ products of the first type and $n-1-s$ products of the second type.
\begin{lemma}
$\Pi_{n,s}$ is upper semimodular.
\end{lemma}
\begin{proof}
Let $a,b,c$ be three elements of $\Pi_{n,s}$ such that $a$ and $b$ cover $c$; denote by 
$M=\{M_1,\ldots,M_t\}$ the set of monomials corresponding to $c$. It is clear that each of the elements $a,b$ is obtained by taking a product of certain \emph{two} monomials. Consider three possible cases.
\begin{enumerate}
\item $a$ is obtained from $c$ by replacing $M_j,M_k$ by $M_j\star_\alpha M_k$, while $b$ is obtained by replacing $M_l,M_m$ by $M_l\star_\beta M_m$, where all indices $j,k,l,m$ are pairwise distinct. Then let $d$ be obtained from $c$ by replacing $M_j,M_k,M_l,M_m$ by $M_j\star_\alpha M_k,M_l\star_\beta M_m$.
\item $a$ is obtained from $c$ by replacing $M_j,M_k$ by $M_j\star_\alpha M_k$, while $b$ is obtained by replacing $M_j,M_l$ by $M_j\star_\beta M_l$ where $j,k,l$ are pairwise distinct. Then let $d$ be obtained from $c$ by replacing $M_j,M_k,M_l$ by $M_j\star_\alpha M_k\star_\beta M_l$.
\item $a$ is obtained from $c$ by replacing $M_j,M_k$ by $M_j\star_\alpha M_k$, while $b$ is obtained by replacing $M_j,M_k$ by $M_j\star_\beta M_k$. Then it is clear that $\alpha\ne\beta$ and that there exists at least one more monomial $M_l$ in $c$ (otherwise $a,b$ would be two distinct maximal elements in the segment, which is a contradiction). Let $d$ be obtained from $c$ by replacing $M_j,M_k,M_l$ by $M_j\star_\alpha M_k\star_\beta M_l$.
\end{enumerate}
In each case $d$ covers both $a$ and $b$.
\end{proof}
\begin{lemma}
Each segment in $\Pi_n$ is isomorphic to the Cartesian product of several maximal segments 
$\Pi_{n_i,r_i}$ (in smaller posets $\Pi_{n_i}$).
\end{lemma}
\begin{proof}
Consider the maximal element of a segment $[x,y]$. Take the corresponding set of monomials and, in particular, the partition $[n]=I_1\sqcup I_2\sqcup\ldots\sqcup I_s$. According to this partition, the segment $[x,y]$ can be decomposed into a product of segments $[x_1,y_1]\times [x_2,y_2]\times\ldots\times[x_s,y_s]$, where the maximal elements $y_t$ are monomials from $\CC(I_t)$.
Let $\tilde{x}=x_k$, $\tilde{y}=y_k$ and denote by $m_1,\ldots,m_r$ the monomials corresponding to $\tilde{x}$. Suppose that the monomial corresponding to $\tilde{y}$ is obtained from these by $s$ products of the first type. Then it is clear that $[\tilde{x},\tilde{y}]$ is isomorphic to the segment $\Pi_{r,s}$ of $\Pi_r$ (replace the letters in the elements of this segment by the corresponding monomials).
\end{proof}

It is clear that the Cartesian product of upper semimodular posets is upper semimodular. Hence $\Pi_n$ is totally upper semimodular. Since any totally upper semimodular poset is Cohen--Macaulay, the theorem follows.
\end{proof}
\begin{corollary}\label{main}\
The operads $\LL$, $\CC$, $\PP$ are Koszul; the $\mathbb{S}$-modules $\PP$ and $\Com\circ\LL$ are isomorphic.
\end{corollary}
\begin{proof}
For the operad $\CC$ its Koszulness is an immediate corollary of our theorem and Proposition \ref{poset-operad}. The operad $\LL$ is Koszul since its dual is Koszul. Both statements about the operad $\PP$, follow from Proposition~\ref{distrib}, since this operad is obtained from Koszul operads $\Com$ and $\LL$ by a distributive law (the Leibniz rule)~--- the proof is similar to the corresponding proof for the Poisson operad in \cite{M}.
\end{proof}

\section{Calculation of the dimensions and characters.}\label{chi-calc}

\subsection{An example of usage of functional equations.}\label{use}
As an example we prove a well known result on the structure of the Lie operad. Since $\L^!=\Com$ and
$$F_\Com(p_1,\ldots,p_n,\ldots)=\sum_{k\ge1}h_k=\exp\left(\sum_{k\ge1}\frac{p_k}{k}\right)-1$$
(here $h_n$ denote the complete symmetric functions (in the case of symmetric polynomials in finitely many variables $h_n$ is a sum of all monomials of degree~$n$), i.e. the Schur polynomials corresponding to the partitions into one part), we have
 $$
\exp\left(-\sum_{k\ge1}\frac{p_k\circ F_\L}{k}\right)-1=-p_1, \mbox{ i.e. }
\sum_{k\ge1}\frac{p_k\circ F_\L}{k}=-\ln(1-p_1).
 $$

Let us prove a formula similar to the classical Moebius inversion formula~\cite{GKP}.

\begin{lemma}[Moebius inversion in the ring $\Lambda$]\label{Moebius}
For $A,B\in \Lambda$ the relation
 $$
A=\sum_{k\ge1}p_k\circ B \quad{\rm and}\quad B=\sum_{k\ge1}\mu_k p_k\circ A
 $$
are equivalent. (Here $\mu$ denotes the number-theoretic Moebius function.)

\end{lemma}

\begin{proof}
Let us obtain, for example, the second relation from the first:
\begin{multline*}
\sum_{k\ge1}\mu_k p_k\circ A=\sum_{k\ge1}\mu_k p_k\circ(\sum_{l\ge1}p_l\circ B)=\\=
\sum_{k\ge1,l\ge1}\mu_k p_{kl}\circ B=\sum_{n\ge1}(\sum_{d\mid n}\mu_d)p_n\circ B=p_1\circ B=p_1.
\end{multline*}
\end{proof}

From this formula we get
$$
F_\L=-\sum_{k\ge1}\frac{\mu_kp_k\circ\ln(1-p_1)}{k}=-\sum_{k\ge1}\frac{\mu_k\ln(1-p_k)}{k}=
\sum_{n\ge1}\frac{1}{n}\sum_{k\mid n}\mu_k p_k^{n/k}.
$$
This relation (proved in another way already in \cite{Br}) is known to be equivalent 
to the following statement~\cite{Kl}.

\begin{theorem}
If the ground field is algebraically closed, the representation of $S_n$ in $\L(n)$ is isomorphic to the representation $\Ind_{H_n}^{S_n}\tau_n$ where $H_n$ denotes a subgroup generated by a cycle of length~$n$, and $\tau_n$ is (any) exact one-dimensional representation of this group.
\end{theorem}

\subsection{Formulas for $SL_2\times S_n$-characters.}

Recall that the components of our operads are equipped by the action of $SL_2$ (arising from the action on the space of generators of the operad~$\LL$), and this action commutes with the symmetric groups. All the information about these operads will follow from the functional equation on the characters and the explicit description (given in Proposition~\ref{CCdim}) of the representation $\CC(n)$ of the group $SL_2\times S_n$.
The character ring of $SL_2$ is isomorphic to the ring of Laurent polynomials in one variable~$q$
(in a way that, for example, the character of the $n$-dimensional irreducible representation is equal to $\frac{q^n-q^{-n}}{q-q^{-1}}$); the element of $\Lambda_{SL_2}$ corresponding to an $\mathbb{S}$-module $\calV$ is denoted by $F_\calV(p_1,\ldots,p_n,\ldots;q)$. This notation differs a little from the one introduced above but we hope that it will not lead to a mess.
In this case the plethysm is defined as follows: $p_n\circ q=q^n$.

Here we list the formulas for the characters of our operads. The proofs are explained in the next section.

Let us introduce the following notation:
\begin{gather*}
a_n(q)=\mu_n\frac{q^n-q^{-n}}{n(q-q^{-1})},\\
c_n(q)=\sum_{d\mid n}\frac{q^d}{d}a_{n/d}(q^d),\\
d_n(q)=\sum_{d\mid n}\frac{1}{d}a_{n/d}(q^d).
\end{gather*}

Now we can formulate our main theorem. We determine the values of $SL_2\times S_n$-characters for the representations in the components of our operads. The value of characters on elements of $S_n$ are $SL_2$-characters, i.e., Laurent polynomials. Take a permutation $\sigma\in S_n$ containing $n_1$ cycles of length~$1$, \ldots, $n_t$ cycles of length~$t$ in the decomposition into disjoint cycles.

\begin{mtheorem}\label{SLchar}
1.~The value of the character of the representation $\PP(n)$ on $\sigma$ is equal to
\begin{multline}
\prod_{s\ge1}\left[\prod_{m=1}^{n_s}\left(\sum_{d\mid s}(sn_dc_{s/d}(q^d))-sm(q^s-q^{-s})+sd_s(q)\right)\right]
\times\\ \times
\prod_{s\ge2}
\frac{\sum\limits_{d\mid s,d\ne s}n_dc_{s/d}(q^d)+d_s(q)}
{\sum\limits_{d\mid s,d\ne s}n_dc_{s/d}(q^d)+n_sq^{-s}+d_s(q)}.
\end{multline}
2.~The value of the character of the space $\LL(n)$ on $\sigma$ can be computed as follows. Let $n_1=n_2=\ldots=n_{k-1}=0$, $n_k\ne0$ (i.e. the length of the shortest cycle is equal to~$k$). Then the character vanishes on $\sigma$ if there exists at least one cycle whose length is not divisible by~$k$; otherwise it is equal to
\begin{multline}
a_k(q)\prod_{j=1}^{n_k-1}\left((n_k-j)q^k+jq^{-k}\right)\times\\
\times\prod_{\substack{s\ge2\colon\\
n_s\ne0}}
\left[\prod_{m=1}^{n_s-1}\left(\sum_{d\colon k\mid d, d\mid s}(sn_{d}c_{s/d}(q^d))-sm(q^s-q^{-s})\right)\times\right.\\ \left.\times
\left(\sum\limits_{d\colon k\mid d, d\mid s, d\ne s}s\,n_{d}c_{s/d}(q^d)\right)
\right].
\end{multline}
\end{mtheorem}

\begin{remark}
Let us compare the character formula for the operad $\PP$ with the character formula for the diagonal harmonic polynomials. From the Haiman's results~\cite{Ha} one can easily deduce that the value of the character of the space of diagonal harmonics in $n$ pairs of variables on a permutation 
(containing $n_1$ cycles of length~$1$, \ldots, $n_k$ cycles of length~$k$ in the decomposition into disjoint cycles) is equal to
\begin{equation}
\frac{(-1)^{\sum_{m\ge1}(m-1)n_m}q^{-n(n-1)/2}}{1+q+q^2+\ldots+q^n}\prod_{m\ge1}(1+q^m+q^{2m}+\ldots+q^{nm})^{n_m}.
\end{equation}
At least, both of the formulas seem rather surprising since the values of the character can be decomposed into small (compared to the values themselves) factors.
\end{remark}

The dimensions, $S_n$-characters and $SL_2$-characters of our operads can be easily obtained from the listed formulas by appropriate specializations. We do not list here the formulas for the $S_n$-characters (corresponding to the specialization $q=1$) since they are approximately of the same complexity. Still the formulas for the dimensions and the $SL_2$-characters are quite simple and beatiful and deserve to be listed.
\begin{corollary}\label{dim}
We have
 $$
f_{\LL}(x)=\sum_{n\ge1}\frac{n^{n-1}x^n}{n!},\qquad f_{\PP}(x)=\sum_{n\ge1}\frac{(n+1)^{n-1}x^n}{n!}.
 $$
Furthermore, denote by $f_\calA(x,q)$ the exponential generating function for the $SL_2$-characters (for an $\mathbb{S}$-module $\calA$ with the $SL_2$ action). Then we have
\begin{gather}
f_{\LL}(x,q)=\sum_{n\ge1}\frac{\prod_{k=1}^{n-1}(kq+(n-k)q^{-1})x^n}{n!},\\
f_{\PP}(x,q)=\sum_{n\ge1}\frac{\prod_{k=1}^{n-1}(kq+1+(n-k)q^{-1})x^n}{n!}.
\end{gather}
\end{corollary}

\begin{proof}
The generating functions for the dimensions are obtained from the characters by setting 
$p_1=x$, $p_k=0$ for $k>1$, $q=1$. The case of $SL_2$-characters is even simpler, since the substitution $q=1$ is not needed.
\end{proof}

\subsection{Proof of the main theorem.}

\subsubsection{The operad~$\LL$.}
Let us note that the $SL_2\times S_n$-character of the operad~$\CC$ is equal to
\begin{multline}
F_{\CC}(p_1,\ldots,p_n,\ldots;q)=\sum_{n\ge1}\frac{q^n-q^{-n}}{q-q^{-1}}h_n=\\
=\frac{\exp\left(\sum_{k\ge1} \frac{q^kp_k}{k}\right)-
\exp\left(\sum_{k\ge1} \frac{q^{-k}p_k}{k}\right)}{q-q^{-1}}.
\end{multline}
To simplify the formulas, let $G:=F_{\LL}(p_1,\ldots,p_n,\ldots;q)$. Theorem~\ref{FunEq} implies that $G$ satisfies the functional equation
 $$
-\frac{1}{q-q^{-1}}\left(\exp\left(\sum_{k\ge1} \frac{-q^kp_k\circ G}{k}\right)-
\exp\left(\sum_{k\ge1} \frac{-q^{-k}p_k\circ G}{k}\right)\right)=p_1.
 $$
Setting $H:=\sum_{k\ge1}\frac{q^k-q^{-k}}{k(q-q^{-1})}p_k\circ G$, we have
 $$
p_1=\frac{\exp((q-q^{-1})H)-1}{q-q^{-1}}\exp\left(-\sum_{k\ge1}\frac{q^kp_k\circ G}{k}\right).
 $$

\begin{lemma}[A $q$-analogue of the Moebius inversion] \label{qMoebuis}
For $A,B\in\Lambda_{SL_2}$ the relations
 $$
A=\sum_{k\ge1}\frac{q^k-q^{-k}}{k(q-q^{-1})}p_k\circ B \quad{\rm and\quad}
B=\sum_{k\ge1}\mu_k\frac{q^k-q^{-k}}{k(q-q^{-1})}p_k\circ A
 $$
are equivalent.
\end{lemma}

\begin{proof}
Analogous to the proof of Lemma~\ref{Moebius}.
\end{proof}

It follows that $\sum_{k\ge1}\frac{q^kp_k\circ G}{k}=\sum_{m\ge1} c_m(q)p_m\circ H$.

Thus,
 $$
p_1=\frac{\exp((q-q^{-1})H)-1}{q-q^{-1}}\exp(-\sum_{m\ge1} c_m(q)p_m\circ H),
 $$
and so
 $$
p_k=p_k\circ p_1=
\frac{\exp((q^k-q^{-k})p_k\circ H)-1}{q^k-q^{-k}}
\exp\left(-\sum_{m\ge1} (c_m(q^k) p_{km}\circ H)\right).
 $$

Let $z_k=p_k\circ H$. Notice that this change of coordinates is triangular and that
\begin{multline*}
\frac{\partial p_k}{\partial z_k}=\exp(z_k(q^k-q^{-k})) \exp(-\sum_{m\ge1} c_m(q^k) z_{km})+\\
+\frac{\exp(z_k(q^k-q^{-k}))-1}{q^k-q^{-k}}(-q^k)
\exp(-\sum_{m\ge1} c_m(q^k) z_{km})=\\
=\exp(-\sum_{m\ge1} c_m(q^k) z_{km})\frac{q^k-q^{-k}\exp(z_k(q^k-q^{-k}))}{q^k-q^{-k}}.
\end{multline*}

We want to compute the coefficient of $p_1^{n_1}\ldots p_k^{n_k}$ in the series $H$. It is equal to the residue (for the sake of briefness, we use the notation $dz=dz_1\ldots dz_k$, $dp=dp_1\ldots dp_k$)
$\res \frac{H(p_1,\ldots,p_n,\ldots;q)}{p_1^{n_1+1}\ldots p_k^{n_k+1}}\,dp$,
and we will now obtain a formula for this residue. Transform it as follows:
\begin{multline*}
\res\prod_{l=1}^k\left(
\exp(-\sum_{m\ge1} c_m(q^l) z_{lm})
\frac{\exp(z_l(q^l-q^{-l}))-1}{q^l-q^{-l}}\right)^{-n_l-1}
z_1\,dp=\\
=\res
\prod_{l=1}^k\left(
\exp(-\sum_{m\ge1} c_m(q^l) z_{lm})
\frac{\exp(z_l(q^l-q^{-l}))-1}{q^l-q^{-l}}\right)^{-n_l-1}
z_1\times\\
\times
\prod_{l=1}^k
\left(\exp(-\sum_{m\ge1} c_m(q^l) z_{lm})
\frac{q^l-q^{-l}\exp(z_l(q^l-q^{-l}))}{q^l-q^{-l}}\right)
\,dz,
\end{multline*}
which is equal to
\begin{multline*}
\res
z_1
\prod_{m=1}^k\left[\exp\left(z_m\sum_{d\mid m}n_dc_{m/d}(q^d)\right)\times\right.\\
\left.\times\left(
\frac{\exp(z_m(q^m-q^{-m}))-1}{q^m-q^{-m}}\right)^{-n_m-1}
\frac{q^m-q^{-m}\exp(z_m(q^m-q^{-m}))}{q^m-q^{-m}}\right]
\,dz.
\end{multline*}

It is obvious that this residue can be decomposed into a product of one-dimensional residues which are computed using Lemma~\ref{resid} (see the Appendix).
Let $x=z_m$, $n=n_m$, $a=\sum_{d\mid m}n_dc_{m/d}(q^d)$,
$b=q^m-q^{-m}$, $\lambda=q^m$, $\mu=q^{-m}$. Consider the cases $m>1$ and $m=1$ separately. In the first case we need to compute the residue
\begin{multline}
\res\frac{\exp(ax)b^{n+1}}{(\exp(bx)-1)^{n+1}}\frac{\lambda-\mu\exp(bx)}{b}\,dx=\\
\res b^n\frac{\lambda\exp(ax)-\mu\exp((a+b)x)}{(\exp(bx)-1)^{n+1}}\,dx,
\end{multline}
which can be (due to Lemma~\ref{resid}) rewritten as
\begin{multline}
b^n\left(\lambda\frac{(a-b)(a-2b)\ldots(a-nb)}{n!b^{n+1}}-
\mu\frac{a(a-b)\ldots(a-(n-1)b)}{n!b^{n+1}}\right)=\\
=\frac{(a-b)(a-2b)\ldots(a-nb)}{n!}\frac{a-n\lambda}{a-nb}.
\end{multline}
(We write the formula in this way to make it correct even for $n=0$~--- then the empty product in the numerator is equal to~$1$ by a standard convention, and the value of the residue is also equal to~$1$, which is true.)

Let now $m=1$. Preserve the same notation. Then the residue is clearly equal to the derivative of the previous residue w.r.t. $a$ at the point 
$\sum_{d\mid 1}n_dc_{m/d}(q^d)=n_1c_1(q)\hm=n_1q\hm=n\lambda$ (according to our notation). Thus the only non-vanishing summand in the derivative corresponds to the factor $(a-n\lambda)$. Hence the residue is equal to $\frac{(a-b)(a-2b)\ldots(a-(n-1)b)}{n!}$ (and zero for  $n=0$).

As a corollary of our calculations we get the following formula for the series $H$:
\begin{multline}\label{LL1}
H(p_1,\ldots,p_n,\ldots;q)=\\=\sum_{k,n_1>0,n_2\ldots n_k\ge0}
\frac{p_1^{n_1}\ldots p_k^{n_k}}{n_1!\ldots n_k!}\prod_{l=1}^{n_1-1}\biggl((n_1-l)q+lq^{-1}\biggr)
\times\\ \times
\prod_{s=2}^k\left[\prod_{l=1}^{n_s}\biggl(\sum_{d\mid s}n_dc_{s/d}(q^d)-l(q^s-q^{-s})\biggr)
\frac{\sum\limits_{d\mid s, d\ne s}n_dc_{s/d}(q^d)}
{\biggl(\sum\limits_{d\mid s,d\ne s}n_dc_{s/d}(q^d)\biggr)+n_sq^{-s}}\right].
\end{multline}

As we mentioned before, the character of the operad~$\LL$ is related to $H$ by the formula
\begin{equation}\label{LL-eq}
F_{\LL}(p_1,\ldots,p_n,\ldots;q)=\sum_{k\ge1} a_k(q)p_k\circ H.
\end{equation}
This character easily implies the formula from the main theorem (to obtain a coefficient of a monomial one does not really need to sum anything).

\subsubsection{The operad~$\PP$.}

Note that we have
\begin{multline*}
F_{\PP}(p_1,\ldots,p_n,\ldots;q)=F_{\Com\circ\LL}(p_1,\ldots,p_n,\ldots;q)=\\=
F_{\Com}(p_1,\ldots,p_n,\ldots)\circ F_{\LL}(p_1,\ldots,p_n,\ldots;q)=\\=
\left(\exp(\sum_{k\ge1} \frac{p_k}{k})-1\right)\circ F_{\LL}(p_1,\ldots,p_n,\ldots;q).
\end{multline*}

We again reduce the calculation of the character values to the calculation of the multidimensional residues. Let $y_i=p_i\circ G$, $z_i=p_i\circ H$ (the notation $G,H$ is introduced above). The change of coordinates $p_i\mapsto y_i$ is uppertriangular, while the change $y_i\mapsto z_i$ is upper unitriangular. Thus to rewrite the residue of the form defined in coordinates $p_1,\ldots,p_k$ in terms of a residue in coordinates $z_1,\ldots,z_k$ we should compute the product $\prod_{l=1}^k \frac{\partial p_l}{\partial y_l}$.
Since $-p_1=\varepsilon(F_{\CC})\circ G$, we have
 $$
-p_k=p_k\circ(\varepsilon(F_{\CC})\circ G)=
F_{\CC}(-y_{k},\ldots,-y_{kn},\ldots;q^k)
 $$
and so $\frac{\partial p_k}{\partial y_k}=
(\frac{\partial\phantom{p_1}}{\partial p_1}F_{\CC})(-y_{k},\ldots,-y_{kn},\ldots;q^k)
$.

It follows that
\begin{multline*}
\res \frac{1}{p_1^{n_1+1}\ldots p_k^{n_k+1}}F_\PP(p_1,\ldots,p_n,\ldots;q)\,dp=\\=
\res\left(\exp\left(\sum_{k\ge1}\frac{y_k}{k}\right)-1\right)
\prod_{l=1}^k\frac{(\frac{\partial\phantom{p_1}}{\partial p_1}F_{\CC})(-y_{l},\ldots,-y_{ln},\ldots;q)}
{(-F_{\CC}(-y_l,\ldots,-y_{ln},\ldots;q))^{n_l+1}}\,dz.
\end{multline*}

In the coordinates $z_i$ the latter residue can be rewritten in the form
\begin{multline*}
\res\left(\exp(\sum_{m\ge1} d_m(q) z_m)-1\right)\times\\ \times
\prod_{l=1}^k\left[\left(\frac{\exp((q^l-q^{-l})z_l)-1}{q^l-q^{-l}}
\exp(-\sum_{m\ge1} c_m(q^l) z_{ml})\right)^{-n_l-1}\right.\times\\
\left.\times
\frac{q^l-q^{-l}\exp((q^l-q^{-l})z_l)}{q^l-q^{-l}}
\exp(-\sum_{m\ge1} c_m(q^l) z_{ml})\right]\,dz.
\end{multline*}

Hence we need to compute the difference of two residues
\begin{multline*}
\res
\prod_{m=1}^k\left[\exp\left(z_m(d_m(q)+\sum_{d\mid m}n_dc_{m/d}(q^d))\right)\times\right.\\
\left.\times\left(
\frac{\exp(z_m(q^m-q^{-m}))-1}{q^m-q^{-m}}\right)^{-n_m-1}
\frac{q^m-q^{-m}\exp(z_m(q^m-q^{-m}))}{q^m-q^{-m}}\right]
\,dz
\end{multline*}
and
\begin{multline*}
\res
\prod_{m=1}^k\left[\exp\left(z_m\sum_{d\mid m}n_dc_{m/d}(q^d)\right)\times\right.\\
\left.\times\left(
\frac{\exp(z_m(q^m-q^{-m}))-1}{q^m-q^{-m}}\right)^{-n_m-1}
\frac{q^m-q^{-m}\exp(z_m(q^m-q^{-m}))}{q^m-q^{-m}}\right]
\,dz,
\end{multline*}
which are similar to the residues from the previous section. It turns out that the second one is equal to zero, and the formula listed above immediately follows.

\section{The monomial basis.}\label{monomial}
Using the dimension formula, we can introduce the basis of the multilinear part of the free algebra over the operad~$\LL$. We use the following result (to appear in~\cite{BDK}).

\begin{proposition}
Given a finite ordered set $A\hm=\{a_1<a_2<\ldots<a_n\}$, define a family of monomials $\mathfrak{B}(A)$ in the free algebra with two compatible brackets generated by $A$ recursively as follows.
\begin{itemize}
\item For $A=\{a_1\}$ let $\mathfrak{B}(A)=\{a_1\}$.
\item If $n>1$ then the monomial $b$ belongs to $\mathfrak{B}(A)$ if and only if it satisfies either of the two conditions:
\begin{enumerate}
\item $b=\{a_i, b'\}_1$, where $i<n$, $b'\in\mathfrak{B}(A\setminus\{a_i\})$;
\item $b=\{b_1,b_2\}_2$, where $b_1\in\mathfrak{B}(A_1)$, $b_2\in\mathfrak{B}(A_2)$
for some $A_1 \sqcup A_2\hm=A$, $a_n\in A_2$, and either $b_1=a_i$ for some~$i$ or 
$b_1$ is a bracket of the first type of two monomials belonging to bases for some subsets of $A_1$.
\end{enumerate}
\end{itemize}
Then $|\mathfrak{B}(A)|=|A|^{|A|-1}$ and the linear span of $\mathfrak{B}(A)$ coincides with the multilinear part of the free algebra.
\end{proposition}

(Firstly, some straightforward arguments show that these monomials span the multilinear part of the free algebra, and then solving some functional equation involving the generating function for the cardinalities of the corresponding sets of monomials gives the statement about the dimensions.)

From the dimension formula we immediately deduce
\begin{theorem}
$\mathfrak{B}(A)$ is a basis of the multilinear part of the free algebra with two compatible brackets generated by $A$.
\end{theorem}

\begin{corollary}
The basis of the multilinear part of the free bihamiltonian algebra generated by $A$ consists of the monomials $b_1\star b_2\star\ldots\star b_k$, where $b_j\in \mathfrak{B}(A_j)$ (for all possible partitions $A=A_1\sqcup \ldots\sqcup A_k$ of $A$; the order of the parts is not important, for example, we can arrange $A_i$ according to the value of the smallest element).
\end{corollary}

\begin{proof}
According to Corollary \ref{main}, $\PP=\Com\circ\LL$, hence the statement follows just from the definition of the composition for $\mathbb{S}$-modules.
\end{proof}

\section{Multiplicities of irreducible representations.}\label{mult}

Here we list some results on the $S_n$-module structure in the multiplicities of irreducible representations of $SL_2$ and the $SL_2$-module structure in the multiplicities of irreducible representations of $S_n$ in the spaces $\LL(n)$ and $\PP(n)$.
We use the following notation: $\mathbbold{1}=\pi_{n}$ and $V=\pi_{n-1,1}$ stand for (respectively) the trivial and the simplicial representation of $S_n$, $L(n)$ denotes the $n$-dimensional irreducible representation of $SL_2$.

\subsection{$S_n$-multiplicities.}

\begin{proposition}
We have the isomorphisms of $SL_2$-modules
\begin{gather}
\label{LLtriv}\Hom_{S_n}(\mathbbold{1},\LL(n))=0,\\
\label{LLsim}\Hom_{S_n}(V,\LL(n))=L(1)^{\otimes (n-1)},\\
\label{PPtriv}\Hom_{S_n}(\mathbbold{1},\PP(n))=L(0),\\
\label{PPsim}\Hom_{S_n}(V,\PP(n))=L(1)\oplus L(1)^{\otimes 2}\oplus\ldots\oplus L(1)^{\otimes (n-1)}.
\end{gather}
\end{proposition}

\begin{proof}
We will discuss in details the case of $\LL$, the case of $\PP$ is completely analogous.
We use the representation in the multilinear part of the free algebra as a realization of the representation $\LL(n)$.

Let us start with the isomorphism \eqref{LLtriv}. 
To isolate the invariants we use the projector $\frac{1}{n!}\sum_{\sigma\in S_n}\sigma$ from the group algebra of the symmetric group. It is enough to prove that this projector annihilates any monomial in the free algebra. Consider an arbitrary monomial $m$; let  $a_i$ and $a_j$ be two generators on the lowest level of bracketing. All permutations in the formula for the projector can be split into pairs 
$(\tau,\tau\cdot(ij))$. It is clear that application of these permutations to $m$ gives two summands differing only by a sign. Hence the projector acts by zero, which is what we want.

Consider now the isomorphism \eqref{LLsim}. Take a subgroup $S_{n-1}$ of $S_n$ preserving $n$. According to the Frobenius reciprocity law, for each representation $T$ of $S_n$ 
$$\Hom_{S_{n-1}}(\mathbbold{1},\res_{S_{n-1}}^{S_n} T)=\Hom_{S_n}(\Ind_{S_{n-1}}^{S_n}\mathbbold{1},T)=
\Hom_{S_n}(\mathbbold{1}\oplus V,T).$$
Applying this to the representation $T=\LL(n)$ with zero multiplicity of the trivial representation, we have
$$\Hom_{S_n}(V,\LL(n))\hm=\Hom_{S_{n-1}}(\mathbbold{1},\res_{S_{n-1}}^{S_n}\LL(n)).$$
Let us find the $S_{n-1}$-invariants in $\LL(n)$. 
We again use the projector from the group algebra to find the invariants. This projector annihilates any monomial that contains two generators $a_i$ and $a_j$ with $i,j<n$ on the lowest level of bracketing. Hence any invariant can contain with nonzero coefficients only left-normed commutators having $a_n$ on the lowest level of bracketing.
All these commutators belong to the basis. For each bracketing all the corresponding monomials should occur with the same coefficients due to invariance. It is easy to see that the $SL_2$-character of the linear span of these elements is equal to the character of the $SL_2$-module $L(1)^{\otimes (n-1)}$, which is what we need.
\end{proof}

\subsection{$SL_2$-multiplicities.}

\begin{proposition}
If the ground field is algebraically closed then we have the isomorphisms of $S_n$-modules (in  notation of Section~\ref{use})
\begin{gather}
\label{LLmax}\Hom_{SL_2}(L(n-1),\LL(n))=\Ind_{H_n}^{S_n}\tau_n,\\
\label{LLnext}\Hom_{SL_2}(L(n-3),\LL(n))=\bigoplus_{k=2}^{n-1}\Ind_{H_k}^{S_n}\tau_k,\\
\label{PPmax}\Hom_{SL_2}(L(n-1),\PP(n))=\Ind_{H_n}^{S_n}\tau_n.
\end{gather}
\end{proposition}

\begin{proof}
Let us use the character formulas. Note that the values of $S_n$-characters of these modules are related in an obvious way to the values of $SL_2\times S_n$-character of the operad $\LL$:
the character of the representation in the space $\Hom_{SL_2}(L(n-1),\LL(n))$ is equal to the coefficient of $q^{n-1}$, and the character of the representation in $\Hom_{SL_2}(L(n-1),\LL(n))$ and in $\Hom_{SL_2}(L(n-3),\LL(n))$ sum up to the coefficient of $q^{n-3}$. It follows from the fact that 
$n-1$ and $n-3$ are, evidently, the maximal weights of the $SL_2$-module $\LL(n)$.

Note that the degree of the Laurent polynomial $c_s(q)$ equals $s$ for $s=1$ and is at most $s-2$ for any other $s$. Thus the degree w.r.t. $q$ of the coefficient of $p_1^{n_1}\ldots p_k^{n_k}$ in the formal power series
$H(p_1,\ldots,p_n,\ldots)$ from~\eqref{LL1} is equal to $n_1-1+\sum_{s>1}(s(n_s-1)+n_s-2k_s)$ for some positive integers $k_s$. The latter sum is equal to $-1+\sum_{s\ge1} sn_s-2\sum_{s>1} k_s=n-1-2\sum_{s>1} k_s$. Hence in the case of $L(n-1)$ we need only the contribution to $F_{\LL(n)}$ of the coefficient of $p_1^n$ in~$H$. In the case of $L(n-3)$ we need also the contribution to $F_{\LL(n)}$ of the coefficient of $p_1^{n-sn_s}p_s^{n_s}$. Let us note that the summands corresponding to the plethysm of $p_i$ and $p_1^{m_1}\ldots p_k^{m_k}$ in the formula for $F_{\LL(n)}$ are of degree $i(-1+\sum_{s\ge1} sm_s-2\sum_{s>1} k_s)+i-1=\sum_{s\ge1}ism_s-2i\sum_{s>1}k_s-1$. This means that the summands which really contribute are those where the plethysm is the plethysm with $p_1$, i.e. the summands from $H$.

To check the formula \eqref{LLmax}, it remains to notice that the vectors of maximal weight in our representation correspond to the monomials where all the brackets are of the same type, and so we can apply Klyachko's theorem.

To prove the formula \eqref{LLnext}, we notice that due to the formula \eqref{LLmax} it remains to prove that the coefficient of $q^{n-3}$ in $F_{\LL(n)}$ is equal to the character of the representation $\bigoplus_{k=2}^{n}\Ind_{H_k}^{S_n}\tau_k$. Since 
$\Ind_{H_k}^{S_n}\tau_k\simeq\Ind_{S_k}^{S_n}\Ind_{H_k}^{S_k}\tau_k$, the latter character is equal (in the ring of the symmetric functions) to $\chi=\sum_{k=2}^n p_1^{n-k}F_{\L(k)}$.
Thus it remains to prove that the correspoding coefficients of our symmetric function coincide with those for $\chi$, which easily follows from the explicit formula~\eqref{LL1}.

The statement about the operad $\PP$ follows immediately, since the highest vectors of weight $n-1$ in $\PP(n)$ belong to the subspace $\LL(n)$.
\end{proof}

\section{Appendix: power series and residues.}\label{basic}

In our calculations we use formal power series and Laurent series, including the series in infinitely many variables. In case of one variable the notation and terminology is standard. As usual, if $f(z)$ is a Laurent series, we call the coefficient of $1/z$ the residue of the formal differential form $f(z)\,dz$ (and denote it by $\res f(z)\,dz$).

The residues which occur in our calculations are of some special type. We compute the general residue of this type:

\begin{lemma}\label{resid}
We have
 $$
\res\frac{\exp(az)}{(\exp(bz)-1)^n}\,dz=\frac{(a-b)(a-2b)\ldots (a-(n-1)b)}{b^n(n-1)!}.
 $$
\end{lemma}

\begin{proof}
Let us note that for $n>1$
 $$
\frac{\exp(az)}{(\exp(bz)-1)^n}\,dz=
\frac{\exp((a-b)z)}{b(1-n)}\cdot d\left(\frac{1}{(\exp(bz)-1)^{n-1}}\right).
 $$
Since the residue of an exact form is equal to zero, we get
 $$
\res\frac{\exp(az)}{(\exp(bz)-1)^n}\,dz=\frac{a-b}{b(n-1)}\res\frac{\exp((a-b)z)}{(\exp(bz)-1)^{n-1}}\,dz,
 $$
and iteration of this relation using the formulas
$\res\frac{\exp(cz)}{\exp(bz)-1}\,dz=\frac{1}{b}$ proves the desired statement.
\end{proof}

For the infinite number of variables formal power series which occur in our calculations belong to the completion of the ring of the symmetric functions. The Laurent polynomials are of the type
$\frac{f(p_1,\ldots,p_n,\ldots)}{(p_1\cdot\ldots\cdot p_k)^N}$ (for positive integers $N,k$), where  $f$ is a formal power series. For a Laurent series $g$ we define the residue of a formal differential form $g(p_1,\ldots,p_n,\ldots)\,dp_1\ldots dp_k$ to be equal to the coefficient of $\frac{1}{p_1\cdot\ldots\cdot p_k}$ in $g$. For the infinite number of variables we use the following version of the residue invariance under the change of coordinates (for which the finite-dimensional case is the same as the infinite-dimensional):
\begin{proposition}
Suppose that the Laurent series $f$ contains only nonnegative powers of $p_i$ with $i>k$. Consider the coordinate change
 $$
p_s\mapsto a_s(t_1,\ldots,t_n,\ldots)=t_s g_s(t_1,\ldots,t_n,\ldots), s\ge1,
 $$
where all $g_s$ are formal power series with~$g_s(0,\ldots,0,\ldots)\ne0$.
Then
 $$
f(p_1,\ldots,p_n,\ldots)\,dp_1\ldots dp_k=f(a_1,\ldots,a_n,\ldots)\det(\tfrac{\partial a_i}{\partial t_j})_{i,j=1,\ldots,k}\,dt_1\ldots dt_k.
 $$
\end{proposition}

\begin{proof}
It is clear that it is enough to prove the equality after the substitutions $p_l=t_l=0$ for all $l>k$. After that we can use the theorems on the residues for the finite-dimensional case.
\end{proof}

{\footnotesize

\medskip

\noindent 
Independent University of Moscow,\\
Bolshoj Vlasievsky per., 11, Moscow, Russia, 119002\\
\texttt{dotsenko@mccme.ru}

\medskip

\noindent
Institute for Theoretical and Experimental Physics,\\
Bolshaya Cheremushkinskaya, 25, Moscow, Russia, 117218\\
\texttt{khorosh@itep.ru}}

\end{document}